\documentclass[a4paper,11pt]{amsart}

       \usepackage{palatino, verbatim}
        \usepackage[latin1]{inputenc}
        \usepackage[T1]{fontenc}
        \usepackage{amsthm, amsfonts}
        \usepackage{amsfonts}
        \usepackage{graphicx}
        \usepackage{amssymb}
        \usepackage{amsmath}
        \usepackage{latexsym}

        \newtheorem{thm}{Theorem}[section]
          \newtheorem{cor}[thm]{Corollary}
          \newtheorem{lem}[thm]{Lemma}
          \newtheorem{prop}[thm]{Proposition}

          \theoremstyle{definition}
          
          \newtheorem{rem}{Remark}
          \newcommand\M{{\mathcal M}}
          \newcommand\im{\mathrm{Im}}
          \newcommand\Pic{\mathrm{Pic }}

\begin{document}

\title
{The Gieseker-Petri divisor in $\M_g$ for $g\leq 13$}
\author{Margherita Lelli--Chiesa}
\address{Humboldt Universit\"at zu Berlin, Institut f\"ur Mathematik, 10099 Berlin}
\email{lelli@math.hu-berlin.de}

\begin{abstract}
The Gieseker-Petri locus $GP_g$ is defined as the locus inside $\M_g$ consisting of curves which violate the Gieseker-Petri Theorem. It is known that $GP_g$ has always some divisorial components and it has been conjectured that it is of pure codimension 1 inside $\M_g$. We prove that this holds true for genus up to 13.
\end{abstract}
\maketitle

\section*{Introduction}
Let $\M_g$ be the coarse moduli space of smooth irreducible projective curves of genus $g$. Given $[C]\in \M_g$ and a line bundle $L$ on $C$, we consider the Petri map
$$
\mu_{0,L}:H^0(C,L)\otimes H^0(C,K_C\otimes L^{-1})\rightarrow H^0(C,K_C).$$
This map has been studied in detail because of its importance in the description of the Brill-Noether varieties $G^r_d(C)$ and $W^r_d(C)$. The most important result in this sense is the Gieseker-Petri Theorem (cf. \cite{gies}, \cite{altra}), which asserts that for the generic curve and for any line bundle on it the Petri map is injective. This implies that if $[C]\in \M_g$ is general and the Brill-Noether number $\rho(g,r,d):=g-(r+1)(g-d+r)$ is nonnegative, then $G^r_d(C)$ is smooth of dimension $\rho(g,r,d)$ and the natural map $G^r_d(C)\rightarrow W^r_d(C)$ is a rational resolution of singularities.
The Gieseker-Petri locus is defined as
$$
GP_g=\{[C]\in\M_g\,\vert\, C\textrm{ does not satisfy the Gieseker-Petri Theorem}\}.
$$
It is conjectured that $GP_g$ has pure codimension $1$ inside $\M_g$; an explanation why this is plausible is given below. The expectation has been proved in genus up to $8$ by Castorena (cf. \cite{casto}, \cite{casto3}). Our main result is:
\begin{thm}\label{thm:principale}
The locus $GP_g$ has pure codimension $1$ inside $\M_g$ for $9\leq g\leq 13$.
\end{thm}
Our strategy is to look at the different components of $GP_g$ determined by the numerical type of the linear series for which the Gieseker-Petri Theorem fails. For values of $g,r,d$ such that  both $r+1$ and $g-d+r$ are at least $2$ we define the Gieseker-Petri locus of type $(r,d)$ as
$$
GP^r_{g,d}:=\{[C]\in\M_g\,\vert\, \exists\textrm{ a base point free }(L,V)\in G^r_d(C)\textrm{ with } \ker\mu_{0,V}\neq 0\},
$$
where $\mu_{0,V}$ denotes the restriction of the Petri map to $V\otimes H^0(C,K_C\otimes L^{-1})$.
Clifford's Theorem, along with Riemann-Roch Theorem, restricts the values of $g,r,d$ for which it is necessary to study the component $GP^r_{g,d}$ to the range $0<2r\leq d\leq g-1$. We also recall that, given $[C]\in GP_g$, at least one of the linear series on $C$ for which the Gieseker-Petri Theorem fails is {\em primitive}, that is, complete and such that both $L$ and $K_C\otimes L^{-1}$ are base point free.

In some cases the codimension of $GP^r_{g,d}$ inside $\M_g$ is known but in general it seems quite difficult to determine the irreducible components of $GP_g$ and control their dimension. When $\rho(g,r,d)<0$, the Petri map corresponding to any $g^r_d$ on a genus $g$ curve cannot be injective for dimension reasons and the study of $GP^r_{g,d}$ essentially coincides with that of the Brill-Noether variety $$\M^r_{g,d}:=\{[C]\in\M_g\,\vert\,W^r_d(C)\neq\emptyset\}.$$
In particular, when $\rho(g,r,d)=-1$, the locus $\M^r_{g,d}$, if nonempty, is an irreducible divisor (cf. \cite{eisi}, \cite{steffen}), known as the Brill-Noether divisor. On the other side, if $\rho(g,r,d)<-1$, the codimension of any component $Z$ of $\M^r_{g,d}$ in $\M_g$ is strictly greater than $1$. If it is true that $GP_g$ has pure codimension $1$ inside $\M_g$, then $Z$ must be contained in some \emph{divisorial}\footnote{By divisorial we will always mean a locus of pure codimension $1$.} component of $GP_g$.

When $\rho(g,r,d)\geq 0$, the Gieseker-Petri locus $GP^r_{g,d}$ can be described as the image, under the natural projection $p:\mathcal{G}^r_d\rightarrow M_g$, of the degeneracy locus $X$ of a map of vector bundles locally defined on $\mathcal{G}^r_d$  and globalizing the Petri map \footnote{Here $\mathcal{G}^r_d$ denotes the stack parametrizing pairs $(C,l)$, where $[C]$ is the isomorphism class of a smooth irreducible projective curve of genus $g$ and $l\in G^r_d(C)$; the map $p$ is the projection on the moduli stack $M_g$.}. Divisoriality of $GP^r_{g,d}$ is suggested by the fact that the expected codimension of $X$ is $\rho(g,r,d)+1$ and it would imply that the restriction of $p$ to $X$ has finite fibers. Farkas proved that $GP^r_{g,d}$ always has a divisorial component if $\rho(g,r,d)\geq 0$ (cf. \cite{gabi1}, \cite{gabi}). However, there are only two cases when $GP^r_{g,d}$ is completely understood. The first one is $GP^1_{g,g-1}$, which can be identified with the locus of curves with a vanishing theta-null and is an irreducible divisor (cf. \cite{tex}). The second case is $GP^1_{g,\frac{g+2}{2}}$, for even genus $g\geq 4$. It has been proved by Eisenbud and Harris (cf. \cite{harri}), that this is a divisor which can be described as the branch locus of the natural map $H_{g,\frac{g+2}{2}}\rightarrow\M_g$ from the Hurwitz scheme $H_{g,\frac{g+2}{2}}$ parametrizing coverings of $\mathbb{P}^1$ of degree $(g+2)/2$ having as source a smooth curve $C$ of genus $g$.

We summarize our results. We show that when $g\leq 13$ the components of $GP_g$ whose codimension is either unknown or strictly greater than $1$ are contained in some divisorial components. Most of the inclusions easily follow from some basic remarks made in the first section. In particular, the components $GP^1_{g,k}$ with $\rho(g,1,k)<-1$ are all contained in the Brill-Noether divisor $\M^1_{g,\frac{g+1}{2}}$ if $g$ is odd, and in the locus $GP^1_{g,\frac{g+2}{2}}$ if $g$ is even.

As a matter of notation, let $\stackrel{\circ}{\M^r_{g,d}}$ be the locus of curves having a primitive $g^r_d$. We define  $$\widetilde{GP}^r_{g,d}:=\{[C]\in\M_g\,\vert\, \exists\textrm{ }(L,V)\in G^r_d(C)\textrm{ with } \ker\mu_{0,V}\neq 0\};$$ notice that here we do not require that $(L,V)$ be base point free. If the Brill-Noether number is either $0$ or $1$, we can prove the inclusion of both $\stackrel{\circ}{\M}^{r+1}_{g,d+1}$ and $\stackrel{\circ}{\M}^{r}_{g,d-1}$  inside $\widetilde{GP}^r_{g,d}$. We use a very recent result, due to Bruno and Sernesi, according to which for values of $g,r,d$ such that $\rho(g,r,d)\geq 0$ and $\rho(g,r+1,d)<0$, the locus $\widetilde{GP}^r_{g,d}$ is divisorial outside its intersection with $\M^{r+1}_{g,d}$ (cf. \cite{sernesi}). As a corollary we obtain that, in even genus, $\widetilde{GP}^1_{g,\frac{g+2}{2}}$ coincides with the divisor $GP^1_{g,\frac{g+2}{2}}$ studied by Eisenbud and Harris. 

In the second paragraph we prove Theorem \ref{thm:principale} in  genera $9,10,11$. In addition to the remarks made in the previous section, we use some well known facts about plane curves. The study of the component $\M^3_{10,9}$ requires extra work: we prove that it is contained in $GP^1_{10,6}$ by remarking that any curve of degree $9$ and genus $10$ in $\mathbb{P}^3$ is either a curve of type $(3,6)$ on a non singular quadric surface or the intersection of two cubic surfaces; linear series on a cubic surface $X$ can be easily written down remembering that $X$ is isomorphic to the blow-up of the projective plane in $6$ points.

In the last paragraph we deal with genera $12$ and $13$. The situation gets more complicated because the methods used before do not enable us to control the codimension of $GP^1_{g,g-2}$. We prove the following theorem:
\begin{thm}\label{thm:brutto}
Let $[C]\in GP^1_{g,g-2}$ be a non hyperelliptic curve with no vanishing theta-null. Let us assume that for any $L\in G^1_{g-2}(C)$ such that $\mu_{0,L}$ is not injective, $L$ is primitive and $K_C\otimes L^{-1}\in W^2_g(C)$ is big. Then $C$ carries only a finite number of $L\in W^1_{g-2}(C)$ for which $\ker\mu_{0,L}\neq 0$.
\end{thm}
This generalizes \cite{casto2}, where it is assumed that the plane model $\Gamma$ of $C$ corresponding to $K_C\otimes L^{-1}$ has only singularities which become nodes after a finite number of blow-ups (in a somewhat  oldfashioned way these are called possibly infinitely near nodes). The idea of our proof is to show that we do not need any assumption on the singularities of $\Gamma$ because the non injectivity of $\mu_{0,L}$ implies that $\Gamma$ has at least one double point, which cannot be a cusp of any order if $[C]\not\in GP^1_{g,g-1}$; then we proceed like in  \cite{casto2}. Theorem \ref{thm:brutto} implies Theorem  \ref{thm:principale} in genus $13$ because no $g^2_{13}$ can be composed with an involution. Instead, for a curve $[C]\in GP^1_{12,10}$ it may happen that a $g^2_{12}$, for which the Petri map is not injective, induces a finite covering of a plane curve of lower genus.  We prove that this can be the case only for $[C]\in GP^1_{12,7}\cup GP^1_{12,8}$ (cf. Theorem \ref{thm:coverings}).

I would like to thank my Ph.D. advisor Gavril Farkas for all the helpful discussions.

\section{Some useful inclusions}
In this section we prove some inclusions among different components of $GP_g$, which enable us to restrict the values of $r$ and $d$ for which the codimension of $GP^r_{g,d}$ must be determined.

We start by stating  the following result, due to Sernesi and Bruno, which exhibits some other divisorial components of $GP_g$:
\begin{thm}\label{thm:hope}
Let $g,r,d$ be integers such that $0<2r\leq d\leq g-1$, $\rho(g,r,d)\geq 0$ and $\rho(g,r+1,d)<0$. Then $\widetilde{GP}^r_{g,d}\setminus (\M^{r+1}_{g,d}\cap \widetilde{GP}^r_{g,d})$ has pure codimension 1 inside $\M_g$.
\end{thm}
The proof of Theorem \ref{thm:hope} has just appeared in \cite{sernesi} and we briefly recall the idea. The condition $\rho(g,r+1,d)<0$ assures that on a generic curve of genus $g$ every $g^r_d$ is complete.  In this situation we consider $\varphi:\mathcal{C}\rightarrow S$ a family of smooth curves of genus $g$ not belonging to $GP^{r+1}_{g,d}$ such that the induced map $S\to \M_g$ is dominant and finite, and the relative scheme $\mathbf{W}^{r,d}_{\mathcal{C}/S}\stackrel{\sigma}{\rightarrow}S$ parametrizing couples $(C_s,L_s)$, with $L_s\in W^r_d(C_s)$ (which in this case implies $h^0(C_s,L_s)=r+1$). The scheme 
$$GP^r_{g,d}(\mathcal{C}/S):=\{s\in S\,\vert\,\varphi^{-1}(s)\in \widetilde{GP}^r_{g,d}\}$$
turns out to be image in $S$ of the degeneracy locus $X_{(r+1)(g-d+r)-1}(\mu)$ of a map of vector bundles $\mu:\mathcal{E}_1\otimes\mathcal{E}_2\rightarrow\mathcal{F}$ defined over $\mathbf{W}^{r,d}_{\mathcal{C}/S}$; if $X_{(r+1)(g-d+r)-1}(\mu)$ is nonempty, then its codimension inside $\mathbf{W}^{r,d}_{\mathcal{C}/S}$ is at most $\rho(g,r,d)+1$. The finiteness of the fibers of the restriction of $\sigma$ to $X_{(r+1)(g-d+r)-1}(\mu)$ follows by a result of Steffen (cf.\cite{steffen}), which can be applied because $\sigma$ is projective and dominant and the sheaf $(\mathcal{E}_1\otimes\mathcal{E}_2)^{\vee}\otimes\mathcal{F}$ is ample relative to $\sigma$, namely it is ample when restricted to any fiber of $\sigma$.

Without the condition $\rho(g,r+1,d)<0$, we could still define the sheaves $\mathcal{E}_1$, $\mathcal{E}_2$ and $\mathcal{F}$ in the same way but $\mathcal{E}_1$ and $\mathcal{E}_2$ would be locally free only when restricted to the open subset $\mathbf{W}^{r,d}_{\mathcal{C}/S}\setminus\mathbf{W}^{r+1,d}_{\mathcal{C}/S}$ . Unfortunately, the restriction of $\sigma$ to $\mathbf{W}^{r,d}_{\mathcal{C}/S}\setminus\mathbf{W}^{r+1,d}_{\mathcal{C}/S}$  is not projective and so Steffen's Theorem cannot be applied in this situation.

We now prove some basic inclusions:
\begin{rem}\label{rem:primo}
For $\rho(g,r-1,d-1)<0$ and $r>1$, we have that:
 $$\M^r_{g,d}\subset \M^{r-1}_{g,d-1}=\widetilde{GP}^{r-1}_{g,d-1}.$$
\end{rem}
\begin{proof}
From any $g^r_d$ we can trivially get a $g^{r-1}_{d-1}$ by subtracting a point outside its base locus.
\end{proof}
Next remark concerns the components $GP^1_{g,k}$:
\begin{rem}\label{rem:secondo}
If $g$ is odd, the following sequence of inclusions holds:
$$\M^1_{g,2}\subseteq \M^1_{g,3}\subseteq\ldots\subseteq\M^1_{g,\frac{g+1}{2}},$$
and $\M^1_{g,\frac{g+1}{2}}$ is a Brill-Noether divisor.

Similarly when $g$ is even we have that:
$$\M^1_{g,2}\subseteq \M^1_{g,3}\subseteq\ldots\subseteq \widetilde{GP}^1_{g,\frac{g+2}{2}}.$$
\end{rem}
\begin{proof}
Cosider $k<\frac{g+1}{2}$ if $g$ is odd and $k<\frac{g+2}{2}$ if $g$ is even. Let $[C]\in\M^1_{g,k}$ and $L$ be a complete $g^1_k$ on $C$. By defining $L':=L\otimes\mathcal{O}_C(P)$ with $P$ a point outside the base locus of $K_C\otimes L^{-1}$, one may prove all the inclusions but $\M^1_{g,\frac{g}{2}}\subset \widetilde{GP}^1_{g,\frac{g+2}{2}}$. When $L$ is a complete $g^1_\frac{g}{2}$ on $C$  with base locus $B$ (not necessarily empty), the Base Point Free Pencil Trick implies both
$$\dim\ker\mu_{0,L}=h^0(C,K_C\otimes L^{-2}\otimes\mathcal{O}_C(B))\geq -\rho(g,1,g/2)=2$$
and
$$\dim\,\ker\mu_{0,L'}=h^0(C,K_C\otimes L^{-2}\otimes\mathcal{O}_C(B-P))\geq 1.$$ 
Thus $L'$ is a $g^1_\frac{g+2}{2}$ on $C$ violating the Gieseker-Petri Theorem and $[C]\in \widetilde{GP}^1_{g,\frac{g+2}{2}}$.
\end{proof}
The following result is a corollary of Theorem \ref{thm:hope}. Together with the previous Remark, it implies that all the loci $GP^1_{g,k}$ such that $\rho(g,1,k)<0$ are  contained in a divisorial component of $GP_g$.
\begin{cor}\label{cor:ultimo}
 In even genus the following equality holds:
$$
\widetilde{GP}^1_{g,\frac{g+2}{2}}=GP^1_{g,\frac{g+2}{2}}.$$
\end{cor}
\begin{proof}
 By Remark \ref{rem:secondo}, we have that $\M^1_{g,\frac{g}{2}}\subset \widetilde{GP}^1_{g,\frac{g+2}{2}}$ and so we can write 
 $$\widetilde{GP}^1_{g,\frac{g+2}{2}}=GP^1_{g,\frac{g+2}{2}}\cup \M^1_{g,\frac{g}{2}},$$
where $GP^1_{g,\frac{g+2}{2}}$ is a divisor on $\M_g$. Furthermore $\M^1_{g,\frac{g}{2}}$ is irreducible and of codimension $2$ in $\M_g$ (cf. \cite{fulton}). Our goal is to show that $\M^1_{g,\frac{g}{2}}\subset GP^1_{g,\frac{g+2}{2}}$.

Theorem \ref{thm:hope} implies that $\widetilde{GP}^1_{g,\frac{g+2}{2}}\setminus\M^2_{g,\frac{g+2}{2}}$ is divisorial, and by Remark \ref{rem:primo} we know that $\M^2_{g,\frac{g+2}{2}}\subset \M^1_{g,\frac{g}{2}}$. It follows that $$\M^1_{g,\frac{g}{2}}\setminus\M^2_{g,\frac{g+2}{2}}\subset GP^1_{g,\frac{g+2}{2}},$$ and the same must be true for its closure. If we show that 
$\M^1_{g,\frac{g}{2}}\setminus\M^2_{g,\frac{g+2}{2}}$ is open in $\M^1_{g,\frac{g}{2}}$, then the irreducibility of $\M^1_{g,\frac{g}{2}}$ implies that $\M^1_{g,\frac{g}{2}}\subset GP^1_{g,\frac{g+2}{2}}$ and we have finished. To end the proof it is enough to remark that the generic curve in $\M^1_{g,\frac{g}{2}}$ has a unique $g^1_\frac{g}{2}$ (cf. \cite{abete2}) while a curve inside $\M^2_{g,\frac{g+2}{2}}$ has at least a $1$-dimensional space of $g^1_\frac{g}{2}$'s (all obtained from a $g^2_\frac{g+2}{2}$ by the subtraction of a point).
\end{proof}
Other useful inclusions come from the following remark:
\begin{rem}\label{rem:terzo}
If $\rho(g,r,d)\in\{0,1\}$, then $\stackrel{\circ}{\M}^{r+1}_{g,d+1}\subset GP^r_{g,d}$ and $\stackrel{\circ}{\M}^r_{g,d-1}\subset\widetilde{GP}^r_{g,d}$.
\end{rem}
\begin{proof}
Assume $\rho(g,r,d)=0$. We fix $[C]\in\stackrel{\circ}{\M}^{r+1}_{g,d+1}$ and $L$ a primitive $g^{r+1}_{d+1}$ on $C$. For any $P\in C$, $L\otimes\mathcal{O}_C(-P)$ is a $g^r_d$ on $C$ and so $G^r_d(C)$ contains $$C':=\{L\otimes\mathcal{O}_C(-P)\,\vert\,P\in C\}\cong C.$$ It follows that $\dim\, T_{L\otimes\mathcal{O}_C(-P)}(G^r_d(C))\geq \dim_{L\otimes\mathcal{O}_C(-P)}G^r_d(C)\geq 1$. By remembering that $$\dim\,T_{L\otimes\mathcal{O}_C(-P)}(G^r_d(C))=\rho(g,r,d)+\dim\,\ker\mu_{0,L\otimes\mathcal{O}_C(-P)}=\dim\,\ker\mu_{0,L\otimes\mathcal{O}_C(-P)},$$ one deduces that $L\otimes\mathcal{O}_C(-P)$ does not satisfy the Gieseker-Petri Theorem. Analogously, given $[C]\in\stackrel{\circ}{\M}^r_{g,d-1}$ and $L$ a primitive, complete $g^r_{d-1}$ on $C$, one defines $$C'':=\{L\otimes\mathcal{O}_C(P)\,\vert\,P\in C\}\cong C$$ and, reasoning as above, proves that $[C]\in \widetilde{GP}^{r}_{g,d}$.

For $\rho(g,r,d)=1$, we consider $[C]\in\stackrel{\circ}{\M}^{r}_{g,d-1}$ and $L$ a primitive $g^{r}_{d-1}$ on $C$. The definition of $C''$ is the same. Since we can assume that $\dim\,G^r_d(C)=1$ (otherwise we could soon conclude that $[C]\in GP^r_{g,d}$), it follows that $C''$ is an irreducible component of $G^r_d(C)$. As $C$ must have a base point free $g^r_d$, there exist components of $G^r_d(C)$ different from $C''$. By the Connectedness Theorem (cf.\cite{alberello}, p. 212), $G^r_d(C)$ is connected. It follows that $G^r_d(C)$ is singular and so $[C]\in \widetilde{GP}^r_{g,d}$. We proceed very similarly if $[C]\in\stackrel{\circ}{\M}^{r+1}_{g,d+1}$.
\end{proof}

\section{Proof of Theorem \ref{thm:principale} in genus $9,10,11$}
In this section we prove that, for genus $g\in\{9,10,11\}$, the Gieseker-Petri locus $GP_g$ is of pure codimension $1$ inside $\M_g$.

Let us fix $g=9$. For $r\in\{4,3\}$ and $2r\leq d\leq 8$ and for $r=2$ and $4\leq d\leq 6$, the Brill-Noether number $\rho(g,r-1,d-1)$ is negative and so, by Remark \ref{rem:primo}, we can restrict our analysis to the components $GP^2_{9,d}$ and $GP^1_{9,k}$ for $d\in\{7,8\}$ and $2\leq k\leq 8$. Moreover, Remark \ref{rem:secondo} implies that $\M^1_{9,k}$ is contained in the Brill-Noether divisor $\M^1_{9,5}$ for $k\leq 4$.

Since $\rho(9,2,7)<0$, we now study $\stackrel{\circ}{\M}^2_{9,7}$. Given $[C]\in\stackrel{\circ}{\M}^2_{9,7}$, if we assume that $C$ does not lie in $\M^1_{9,5}$, then any $g^2_7$ on $C$ is base point free and defines an embedding $$\phi:C\rightarrow\Gamma\subset\mathbb{P}^2,$$ where $\Gamma$ is a plane curve of degree $7$ and genus $9$. By the Genus Formula it follows that $\Gamma$ is singular, which is a contradiction.

Regarding the component $GP^2_{9,8}$,  we note that $\rho(9,2,8)=0$ and $\rho(9,3,8)<0$, so Theorem \ref{thm:hope} implies that $\widetilde{GP}^2_{9,8}\setminus (\M^3_{9,8}\cap \widetilde{GP}^2_{9,8})$ is divisorial. We do not need to study $\M^3_{9,8}\cap GP^2_{9,8}$ separately because, by Remark \ref{rem:primo}, the inclusion $\M^3_{9,8}\subseteq \M^2_{9,7}$ holds.

Let us consider the components $GP^1_{9,k}$ for $k\in\{6,7,8\}$. For $k\in\{6,7\}$ we have that $\rho(9,1,k)>0$ and $\rho(g,2,k)<0$ and so the locus $\widetilde{GP}^1_{g,k}\setminus (\widetilde{GP}^1_{g,k}\cap\M^2_{g,k})$ is divisorial. As $GP^1_{9,8}$ is the irreducible divisor consisting of curves with a vanishing theta-null, Theorem \ref{thm:principale} is proved in genus $9$.\vspace{0.7cm}

Before dealing with the case of genus $10$, we prefer to treat the case of genus $11$, which is very similar to the one we have just studied. As before, by applying Remark \ref{rem:primo} and Remark \ref{rem:secondo} we reduce to considering the components $GP^2_{11,d}$ and $GP^1_{11,k}$ for $8\leq d\leq 10$ and $7\leq k\leq 10$.

We can prove that $\stackrel{\circ}{\M}^2_{11,8}$ is contained in the Brill-Noether divisor $\M^1_{11,6}$ simply by remarking that any $g^2_8$ on a genus 11 curve $[C]\not\in\M^1_{11,6}$ is base point free and defines an embedding $$\phi:C\rightarrow \Gamma\subset\mathbb{P}^2.$$ We get a contradiction because $\Gamma$ is a plane curve of degree $8$ and genus $11$ and so it must be singular by the Genus Formula.

Concerning the other components, the locus $\M^2_{11,9}$ is a Brill-Noether divisor, while $\widetilde{GP}^2_{11,10}$ is divisorial outside its intersection with $\M^3_{11,10}$ because $\rho(11,2,10)>0$ and $\rho(11,3,10)<0$. 

Theorem \ref{thm:hope} can be applied in order to prove that the locus $\widetilde{GP}^1_{11,k}\setminus (\M^2_{11,k}\cap \widetilde{GP} ^1_{11,k})$ is divisorial for $7\leq k\leq 9$, too. The component $GP^1_{11,10}$ is the irreducible divisor of curves with a vanishing theta-null and so Theorem \ref{thm:principale} is proved in genus $11$.\vspace{0.7cm}

We now deal with the case of genus $10$. As above, by Remarks \ref{rem:primo} and \ref{rem:secondo}, the only components of $GP_{10}$ we have to consider are $GP^2_{10,d}$ and $GP^1_{10,k}$ for $7\leq d\leq 9$ and $7\leq k\leq 9$.

As $\rho(10,1,6)=0$, Remark \ref{rem:terzo} implies that $\stackrel{\circ}{\M}^2_{10,7}\subset GP^1_{10,6}$. Moreover, $\rho(10,2,9)=1$ and so Remark \ref{rem:terzo} implies that $\stackrel{\circ}{\M}^2_{10,8}\subset\widetilde{GP}^2_{10,9}$, too. Since $\rho(10,3,9)<0$, the locus  $\widetilde{GP}^2_{10,9}$ is divisorial outside $\M^3_{10,9}$. In this case we have to study the component $\M^3_{10,9}$ separately because our remarks imply only that $\stackrel{\circ}{\M}^3_{10,9}\subseteq \stackrel{\circ}{\M}^2_{10,8}\subseteq GP^2_{10,9}$. We postpone the study of $\stackrel{\circ}{\M}^3_{10,9}$. For $k\in\{7,8\}$, the locus $\widetilde{GP}^1_{10,k}\setminus(\widetilde{GP}^1_{10,k}\cap\M^2_{10,k})$ is divisorial because $\rho(10,2,k)<0$, while $GP^1_{10,9}$ is the irreducible divisor consisting of curves with a vanishing theta-null.

In order to end the proof of Theorem \ref{thm:principale} in genus $10$, we now study $\M^3_{10,9}$. We consider $[C]\in\M^3_{10,9}$ and $L$ a $g^3_9$ on $C$. We can assume $[C]\not\in\M^3_{10,8}$ and so $L$, being base point free, defines a morphism $\phi:C\rightarrow\Gamma\subset\mathbb{P}^3$. Furthermore, we can assume that $[C]\not\in\M^2_{10,7}$, which forces $\phi$ to be an embedding. Therefore $C$ can be seen as a curve of genus $10$ and degree $9$ in $\mathbb{P}^3$. By the classification of curves in $\mathbb{P}^3$, we know that $C$ is either a curve of type $(3,6)$ on a non singular quadric surface $S$ or the intersection of two cubic surfaces (cf. \cite{harti2} Example 6.4.3. chp.IV). In the first case the lines of type $(0,1)$ on $S$ cut out a $g^1_3$ on $\Gamma$. The second case is treated in the following lemma:
\begin{lem}
Let $[C]\in\M_{10}$ be the intersection of two cubic surfaces $X,Y$ in $\mathbb{P}^3$. Then $[C]\in GP^1_{10,6}$.
\end{lem}
\begin{proof}
It is classically known that $X$ is isomorphic to the blow-up of $\mathbb{P}^2$ in $6$ points $P_1,\ldots,P_6$. We denote by $\pi:X\rightarrow \mathbb{P}^2$ the projection and by $E_i$ the exceptional divisors. $\Pic(X)\cong\mathbb{Z}^7$ and it is generated by $l,e_1,e_2,\ldots,e_6$, where $l$ is the class of the strict transform of a line in $\mathbb{P}^2$ and $e_i$ is the class of $E_i$. The class of the hyperplane section is $h=3l-\sum e_i$, while $$K_X\sim -h=-3l+\sum e_i.$$ As $C$ lies on another cubic surface $Y$, then $$C\sim 3h=9l-3\sum e_i,$$ namely $C$ is the strict transform of a plane curve $\widetilde{C}$ of degree $9$ with $6$ triple points. The pencil of cubics through $P_1,\ldots,P_6$ with a double point in $P_1$ cuts out a $g^1_6$ on $\widetilde{C}$. The strict transforms of these cubics cut out on $C$ the linear series $$L:=(3l-\sum_{i\neq 1}e_i-2e_1)\otimes\mathcal{O}_C.$$
In order to check that $L$ is a $g^1_6$ on $C$, we tensor with $\mathcal{O}_X(3l-\sum_{i\neq 1}e_i-2e_1)$ the exact sequence
$$0\rightarrow\mathcal{O}_X(-C)\rightarrow\mathcal{O}_X\rightarrow\mathcal{O}_C\rightarrow 0,$$
getting
$$0\rightarrow\mathcal{O}_X(-6l+2\sum_{i\neq 1}e_i+e_1)\rightarrow\mathcal{O}_X(3l-\sum_{i\neq 1}e_i-2e_1)\rightarrow \mathcal{O}_C(3l-\sum_{i\neq 1}e_i-2e_1)\rightarrow 0.$$
As $6l-2\sum_{i\neq 1}e_i-e_1$ is ample (cf. \cite{harti2} Cor.4.13 chap.V), Kodaira Vanishing Theorem implies that $h^i(X,\mathcal{O}_X(-6l+2\sum_{i\neq 1}e_i+e_1))=0$ for $i=0,1$. It follows that
$$
\begin{array}{llll}
h^0(C,\mathcal{O}_C(3l-\sum_{i\neq 1}e_i-2e_1))&=&h^0(X,\mathcal{O}_X(3l-\sum_{i\neq 1}e_i-2e_1))&=\\
&=&h^0(\mathbb{P}^2,\mathcal{O}_{\mathbb{P}^2}(3)\otimes\mathcal{O}_{\mathbb{P}^2}(-\sum_{i\neq 1}P_i-2P_1))&=\\
&=&2&
\end{array}
$$
and this is enough to conclude that $L$ is a pencil on $C$; it is trivial to check that its degree is $6$.\\
By the Base Point Free Pencil Trick we have that
$$
\ker\mu_{0,L}\cong H^0(C,K_C\otimes L^{-2})=H^0(C,\mathcal{O}_C(2e_1)).
$$
As $\mathcal{O}_C(2e_1)$ is effective, it follows that $[C]\in GP^1_{10,6}$.
\end{proof}
\begin{rem}
 The previous Lemma can also be proved by using the results of \cite{marti}. Curves of genus $10$ which are the complete intersection of two cubic surfaces in $\mathbb{P}^3$ are the only curves of Clifford dimension $3$. Martens proved that such curves are $6$-gonal and carry a one-dimensional family of $g^1_6$. Since $\rho(10,1,6)=0$, this is enough to conclude that they lie in $GP^1_{10,6}$.

It is natural to ask whether all curves of Clifford dimension greater than $1$ lie in a divisorial component of the Gieseker-Petri locus. Curves of Clifford dimension $2$ are smooth plane curves of degree $d\geq 5$. Their gonality is $d-1$ and  there is a one-dimensional family of pencils computing it. As $\rho\left(\binom{d-1}{2},1,d-1\right)\leq 0$ for $d\geq 5$, Remark \ref{rem:secondo} implies that they lie in the Brill-Noether divisor $\M^1_{g,\frac{g+1}{2}}$ when $g= \binom{d-1}{2}$ is odd, and in the irreducible divisor $GP^1_{g,\frac{g+2}{2}}$ when $g$ is even.

It is conjectured in \cite{cliff} that if $C$ is a curve of Clifford dimension $r>3$, then $g(C)=4r-2$, $\mathrm{gon}(C)=2r$ and there is a one-dimensional family of pencils computing the gonality (this conjecture was proved in \cite{cliff} for $r\leq 9$). Since $\rho(4r-2,1,2r)=0$, such curves lie in the divisor $GP^1_{g,\frac{g+2}{2}}=GP^1_{4r-2,2r}$.
\end{rem}

\section{Proof of Theorem \ref{thm:principale} in genus $12,13$}
The situation in genus $12$ and $13$ is slightly more complicated as there is a component in $GP_g$ which cannot be studied using the methods explained in the previous sections.

In genus $12$, by Remarks \ref{rem:primo} and \ref{rem:secondo}, we have to analyze only the components $\stackrel{\circ}{\M}^3_{12,11}$, $GP^2_{12,d}$ for $8\leq d\leq 11$ and $GP^1_{12,k}$ for $8\leq k\leq 11$. Since $\rho(12,2,10)=0$, Remark \ref{rem:terzo} implies that both $\stackrel{\circ}{\M}^3_{12,11}$ and $\stackrel{\circ}{\M}^2_{12,9}$ are contained in $\widetilde{GP}^2_{12,10}$. Remark \ref{rem:terzo} can also be used in order to show that $\stackrel{\circ}{\M}^2_{12,8}\subset GP^1_{12,7}$; indeed, $\rho(12,1,7)=0$.

As $\rho(12,3,d)<0$ for $d\in\{10,11\}$, the loci $\widetilde{GP}^2_{12,10}$ and $\widetilde{GP}^2_{12,11}$ are divisorial outside their intersection with $\M^3_{12,10}$ and $\M^3_{12,11}$ respectively. We have to study $\M^3_{12,10}$ separately because our remarks only imply that $\stackrel{\circ}{\M}^3_{12,10}\subset \stackrel{\circ}{\M}^2_{12,9}\subset GP^2_{12,10}$.\\
Given $[C]\in\M^3_{12,10}$, we can suppose that $[C]\not\in \M^2_{12,8}$ and so any $g^3_{10}$ on $C$ is base point free and defines an embedding $\phi:C\rightarrow \Gamma\subset\mathbb{P}^3$. It can be seen that $\Gamma$ has ten $4$-secant lines (cf. \cite{alberello}, p. 351), each of which corresponds to a $g^1_6$ on it.

Theorem \ref{thm:hope} can be applied in order to show that the locus $\widetilde{GP}^1_{12,k}$ is divisorial outside $\M^2_{12,k}$ for $k\in\{8,9\}$. The component $GP^1_{12,11}$ is an irreducible divisor. We postpone the study of $GP^1_{12,10}$ to the end of the section.\vspace{0.7cm}

The situation in genus $13$ is very similar to that in genus $12$. By Remarks \ref{rem:primo} and \ref{rem:secondo}, we reduce to considering $\stackrel{\circ}{\M}^3_{13,12}$, $GP^2_{13,d}$ for $9\leq d\leq 12$ and $GP^1_{13,k}$ for $8\leq k\leq 12$.

As $\rho(13,2,11)=1$, Remark \ref{rem:terzo} implies that both $\stackrel{\circ}{\M}^3_{13,12}$ and $\stackrel{\circ}{\M}^2_{13,10}$ are contained in $\widetilde{GP}^2_{13,11}$.

Concerning $\M^2_{13,9}$, any $g^2_9$ on a genus $13$ curve $[C]\not\in \M^1_{13,7}$ defines an embedding $\phi:C\rightarrow \Gamma\subset\mathbb{P}^2$. We get a contradiction because the Genus Formula forces $\Gamma$ to be singular.

The components $\widetilde{GP}^2_{13,11}$ and $\widetilde{GP}^2_{13,12}$ are divisorial outside $\M^3_{13,11}$ and $\M^3_{13,12}$ respectively. As before we have to study $\M^3_{13,11}$ separately. Given $[C]\in \M^3_{13,11}$ such that $[C]\not\in\M^2_{13,9}$, by taking the $4$-secant lines to the space model of $C$ corresponding to any $l\in G^3_{11}(C)$, one shows that $C$ has a $g^1_7$.

Regarding the other components, the locus $\widetilde{GP}^1_{13,k}$ is divisorial outside its intersection with $\M^2_{13,k}$ for $k\in\{8,9,10\}$, while $GP^1_{13,12}$ is an irreducible divisor. Therefore Theorem \ref{thm:principale} is proved also in genus $13$ if we are able to verify that the component $GP^1_{g,g-2}$ is divisorial. In order to show this, we generalize a result of Castorena (cf. \cite{casto2}) as follows.

We consider curves $[C]\in GP^1_{g,g-2}$ such that for any $L\in G^1_{g-2}(C)$ with $\ker\mu_{0,L}\neq 0$ the following are satisfied:
\begin{enumerate}
\item $L$ is primitive.
\item The morphism $\phi:=\phi_{K_C\otimes L^{-1}}$ is birational.
\end{enumerate}
We remark that the first condition is satisfied if $[C]\not\in GP^1_{g,g-3}\cup GP^2_{g,g-2}\cup GP^2_{g,g-1}$, because if $L$ were not complete (respectively not base point free), this would imply $[C]\in GP^2_{g,g-2}$ (resp. $[C]\in GP^1_{g,g-3}$). Similarly, if $K_C\otimes L^{-1}$ is not base base point free, then $[C]\in GP^2_{g,g-1}$.
We prove the following result:
\begin{prop}\label{prop:bo}
Let us consider $Z_g\subset GP^1_{g,g-2}$ the locus consisting of curves $[C]\in GP^1_{g,g-2}$ such that if $L\in G^1_{g-2}(C)$ satisfies $ker\mu_{0,L}\neq 0$, then $L$ is primitive and $K_C\otimes L^{-1}$ is big. The scheme $Z_g$ has pure codimension $1$ in $\M_g$ outside its intersection with the hyperelliptic locus and the divisor $GP^1_{g,g-1}$.
\end{prop}
It is clear that $Z_g$ is open in $GP^1_{g,g-2}$. In order to prove Proposition \ref{prop:bo} we need the following: 
\begin{lem}\label{lem:coffee}
If $[C]\in Z_g$ and $L$ is a $g^1_{g-2}$ on $C$ such that $\ker\mu_{0,L}\neq 0$, then $L$ is the pullback to $C$ of the pencil cut out on $\Gamma:=\phi_{K_C\otimes L^{-1}}(C)$ by the lines through a singular point $x$. In particular, $x$ is a double point of $\Gamma$ and $K_C\otimes L^{-2}=\frac{1}{k}\phi^*\mathcal{O}_\Gamma(x)$, where $k$ is the number of blow-ups necessary to desingularize $\Gamma$ in $x$ (e.g., if $x$ is a tacnode, then $k=2$). 
\end{lem}
\begin{proof}
The Base Point Free Pencil Trick implies that $K_C\otimes L^{-2}=\phi^*(\mathcal{O}_{\mathbb{P}^2}(1))\otimes L^{-1}$ is effective. Hence $L$ is the pullback to $C$ of the pencil cut out on $\Gamma$ by the lines through a singular point $x$, which must be a double point because $L$ is base point free. Furthermore, $K_C\otimes L^{-2}$ is linearly equivalent to $\frac{1}{k}\phi^*\mathcal{O}_\Gamma(x)$.
\end{proof}
We can now prove the following fact:
\begin{lem}\label{lem:forse}
If $[C]\in Z_g$, $[C]\not\in GP^1_{g,g-1}$ and $C$ is not hyperelliptic, then there exists only a finite number of $L\in W^1_{g-2}(C)$ such that $\mu_{0,L}$ is not injective.
\end{lem}
\begin{proof}
We recall and adapt the proof of Castorena, referring to \cite{casto2} for further details.
Given $L$ a $g^1_{g-2}$ on $C$ with $\ker\mu_{0,L}\neq 0$, we have that $$K_C\otimes L^{-2}=\frac{1}{k}\phi^*\mathcal{O}_\Gamma(x)=\mathcal{O}_C(P+Q),$$ and we can assume $P\neq Q$ because otherwise $L\otimes\mathcal{O}_C(P)$ would be a theta characteristic with a 2-dimensional space of sections, thus contradicting  $[C]\not\in GP^1_{g,g-1}$. We remark that asking that $P\neq Q$ is equivalent to requiring that $x$ be not a cusp of any order. As $C$ is not hyperelliptic, $h^0(C,\mathcal{O}_C(P+Q))=1$ and $h^0(C,K_C\otimes\mathcal{O}_C(-P-Q))=g-2$. It follows that $L^{2}$ lies in the intersection of the following two subvarieties of $\mathrm{Pic}^{2g-4}(C)$:
$$
\begin{array}{lll}
X_1&:=&\{L^{2}\,\vert\,L\in W^1_{g-2}(C)\},\\
X_2&:=&\{K_C\otimes\mathcal{O}_C(-P-Q)\,\vert\, P,Q\in C\}\subset W^{g-3}_{2g-4}(C).
\end{array}
$$
In order to show that $X_1\cap X_2$ is a finite set, it is enough to prove that the intersection 
$T_{L^{2}}(X_1)\cap T_{L^{2}}(X_2)=\{0\}$ in $H^1(C,\mathcal{O}_C)=T_{L^{2}}(\mathrm{Pic}^{2g-4}(C))$, or equivalently, that $T_{L^{2}}(X_1)^\perp +T_{L^{2}}(X_2)^\perp$ generates the whole $H^0(C,K_C)=T_{L^{2}}(\mathrm{Pic}^{2g-4}(C))^\perp$.
Indeed, it is trivial to see that
$$\dim \,T_{L^{2}}(X_1)=\dim \,T_L(W^1_{g-2}(C))=\rho(g,1,g-2)+\dim\,\ker\mu_{0,L}=g-5$$
while $\dim\,T_{L^{2}}(X_2)=2$, because $\mu_{0,L^{2}}$ is injective.\\ We recall that $T_{L^{2}}(X_1)^\perp\simeq \im\,\mu_{0,L}$ and
$$T_{L^{2}}(X_2)^\perp\simeq \im\,\mu_{0,L^{2}}\simeq H^0(C,K_C\otimes\mathcal{O}_C(-P-Q)).$$ We should prove that $\dim\,T_{L^{2}}(X_1)^\perp\cap T_{L^{2}(}X_2)^\perp=3$, that is, $P$ and $Q$ impose independent conditions on the $5-$dimensional space $\mathcal{L}:=\im\,\mu_{0,L}$. As explained in \cite{casto2}, it is enough to show that $P$ and $Q$ impose independent conditions on $D+\vert L\vert\subset\mathcal{L}$, where $D\in \vert K_C\otimes L^{-1}\vert$ is a divisor not containing $P+Q$. 
Indeed, let us consider the finite sequence of blow-ups
$$
X_k\rightarrow X_{k-1}\rightarrow \ldots\rightarrow X_0=\mathbb{P}^2
$$ 
necessary to desingularize $\Gamma$ in $x$. We denote by $C_h$ the strict transform of $\Gamma$ in $X_h$ and by $\phi_h:X_h\rightarrow \mathbb{P}^2$ the projection; we note that $C$ coincides with $C_k$ and the normalization map $\phi_k\vert_{C_k}$ is $\phi_{K_C\otimes L^{-1}}$. The curve $C_{k-1}$ has a node in the point $x_{k-1}$ which maps to $x$ via $\phi_{k-1}\vert_{C_{k-1}}$.  The strict transform of a line through $x$ intersects $C$ along a divisor of the form $E_l\otimes\mathcal{O}_C(P+Q)$, with $E_l\in \vert L\vert$. As $C_{k-1}$ has a node in $x_{k-1}$, there exist two lines $l_1$ and $l_2$ through $x$ whose strict transforms in $X_{k-1}$ are the two tangent lines to $C_{k-1}$ in $x_{k-1}$. It follows that the strict transforms of $l_1$ and $l_2$ in $X_k$ intersect $C$ in $E_1\otimes\mathcal{O}_C(2P+Q)$ and $E_2\otimes\mathcal{O}_C(P+2Q)$ respectively, where $E_1\otimes\mathcal{O}_C(P)\in\vert L\vert$ does not contain $Q$ and $E_2\otimes\mathcal{O}_C(Q)\in \vert L\vert$ does not contain $P$. It follows that $P$ and $Q$ impose independent conditions on $D+\vert L\vert$.
\end{proof}
\begin{proof}[Proof of Proposition \ref{prop:bo}]
Let $[C]\in Z_g$ be a non hyperelliptic curve with no vanishing theta-null. One may find a neighborhood $U\subset\M_g$ of $C$, intersecting neither the hyperelliptic locus nor the divisor $GP^1_{g,g-1}$, such that there exists a finite ramified covering $\pi:\widetilde{U}\rightarrow U$, a universal curve $\varphi:\Gamma_{\widetilde{U}}\rightarrow \widetilde{U}$ and a variety $\mathcal{G}^1_{g-2}\stackrel{\xi}{\rightarrow}\widetilde{U}$ proper over $\widetilde{U}$ which parametrizes pairs $(C,(V,L))$ with $[C]\in\widetilde{U}$ and $(V,L)$ a $g^1_{g-2}$ on $\varphi^{-1}(C)$. Up to restricting $U$, we can also assume that $U\cap GP^1_{g,g-2}\subset Z_g$. The scheme $\mathcal{G}^1_{g-2}$ is smooth of dimension $\rho(g,1,g-2)+\dim\,\M_g$ (cf. \cite{abete}). We define the following subvariety of $\mathcal{G}^1_{g-2}$:
$$
\widetilde{Z}_g:=\{(C,L)\in\mathcal{G}^1_{g-2}\,\vert\,[C]\in\pi^{-1}(Z_g\cap U), \ker\mu_{0,L}\neq 0\}.
$$
Lemma \ref{lem:forse} implies that the fiber of the projection from $\widetilde{Z}_g$ on $Z_g\cap U$ is finite. For any $(C,L)\in\widetilde{Z}_g$, the curve $C$ is not hyperelliptic and so $\dim\im\mu_{0,L}=5$. Locally the Petri map defines a homomorphism $\mu$ of vector bundles on $\mathcal{G}^1_{g-2}$ and $\widetilde{Z}_g$ can be identified with the fifth degeneracy locus of $\mu$. By the fact that each irreducible component of $\widetilde{Z}_g$ has codimension $\leq\rho(g,1,g-2)+1$ in $\mathcal{G}^1_{g-2}$ and by the finiteness of the fibers of $\pi\circ\xi$ over the points of $\pi\circ\xi(\widetilde{Z}_g)$, we can deduce that each component of $Z_g\cap U$ has codimension at most $1$ in $U$. It must be $1$ because of the Gieseker-Petri Theorem.
\end{proof}
As a consequence we gain the following:
\begin{cor}
The locus $GP_{13}$ has pure codimension $1$ in $\M_{13}$.
\end{cor} 
\begin{proof}
By the above discussion we should only study the component $GP^1_{13,11}$. Given $[C]\in GP^1_{13,11}$, we assume that $[C]$ does not lie in $GP^{1}_{13,10}\cup GP^2_{13,11}\cup GP^2_{13,12}$. In particular, condition $(1)$ is satisfied for any $L\in G^1_{13}(C)$ for which the Gieseker-Petri Theorem fails. Moreover,  $K_C\otimes L^{-1}$ cannot be composed with any involution and so condition $(2)$ is satisfied, too. It follows that $[C]\in Z_g$ and so Proposition \ref{prop:bo} is enough to conclude.
\end{proof}
Next we turn to the case of genus $12$. Given $[C]\in GP^1_{12,10}$ such that condition $(1)$ is satisfied for any $L\in G^1_{10}(C)$ with $\ker\mu_{0,L}\neq 0$, it could still happen that some of the above  $L\in W^1_{10}(C)$ violate condition $(2)$, that is, $K_C\otimes L^{-1}$ is not big. We prove the following:
\begin{thm}\label{thm:coverings}
Let $[C]\in GP^1_{12,10}$ and let us assume that condition $(1)$ is satisfied for any $L\in G^1_{10}(C)$ such that $\ker\mu_{0,L}\neq 0$. If for one of such $L\in W^1_{10}(C)$ the morphism $K_C\otimes L^{-1}$ defines a finite covering of a plane curve $\Gamma$ of degree strictly less than $12$, then $[C]$ lies in $GP^1_{12,7}\cup GP^1_{12,8}$.
\end{thm}
\begin{proof}
Let $[C]\in GP^1_{12,10}$ be as in the hypothesis and $L$ be a $g^1_{10}$ on $C$ for which the Gieseker-Petri Theorem fails. If $\phi:=\phi_{K_C\otimes L^{-1}}:C\to \Gamma\subset\mathbb{P}^2$ is not birational, then it is a finite covering of degree $6$, $4$, $3$ or $2$. We analyze these cases.\vspace{0.3cm}

\textbf{(I):} $\deg\phi_{K_C\otimes L^{-1}}=6$. In this case $\Gamma$ is rational and so $C$ has a $g^1_6$.\vspace{0.3cm}

\textbf{(II):} $\deg\phi_{K_C\otimes L^{-1}}=3$. Then $\Gamma$ has degree $4$ and genus at most $3$. If $g(\Gamma)<3$, then $\Gamma$ has at least one singular point and by taking the lines through it one sees that $\Gamma$ is hyperelliptic and so $C$ has a $g^1_6$.

Let us consider $g(\Gamma)=3$. As the triple cover is induced by $K_C\otimes L^{-1}$, it follows that
$
K_C\otimes L^{-1}=\phi^*\mathcal{O}_\Gamma(1)=\phi^*K_\Gamma$ and so $L=\mathcal{O}_C(R)$, where $R$ is the ramification locus. The Base Point Free Pencil Trick thus implies that  
$$
\ker\mu_{0,L}\simeq H^0(C,K_C\otimes\mathcal{O}_C(-2R))\simeq H^0(C,\phi^*\mathcal{O}_\Gamma(1)\otimes\mathcal{O}_C(-R)).
$$
If this were not zero, then there would exist a divisor $D$ on $\Gamma$, $\mathcal{O}_\Gamma(D)=\mathcal{O}_\Gamma(1)$, such that $\phi^*D-R\geq 0$. This would imply that $D$ contains the branch locus $B$ but this is impossible because $\deg B\geq \frac{1}{2}\deg R=5$ while $\deg D=4$. \vspace{0.3cm}

\textbf{(III):} $\deg\phi_{K_C\otimes L^{-1}}=4$. The curve $\Gamma$ has degree $3$ and so it is either a rational curve or a smooth elliptic curve. In the first case $C$ has a $g^1_4$.

If $\Gamma$ is elliptic , then we have that $K_C\otimes L^{-1}=\phi^*\mathcal{O}_\Gamma(1)$ and $$L=\phi^*(K_\Gamma\otimes\mathcal{O}_\Gamma(-1))\otimes\mathcal{O}_C(R)=\phi^*\mathcal{O}_\Gamma(-1)\otimes\mathcal{O}_C(R).$$ It follows that
$$\ker\mu_{0,L}\simeq H^0(C,\mathcal{O}_C(R)\otimes(\mathcal{O}_C(R)\otimes\phi^*\mathcal{O}_\Gamma(-1))^{-2})=H^0(\phi^*\mathcal{O}_\Gamma(2)\otimes\mathcal{O}_C(-R)).$$
This is nonzero whenever there exists a divisor $D$ on $\Gamma$ such that $\mathcal{O}(D)=\mathcal{O}_\Gamma(2)$ and $\phi^*D-R\geq 0$. This never happens because $D$ has degree $6$ and it should contain the base locus $B$, whose degree is at least $\frac{1}{3}\deg R>7$. It follows that there exixts no $g^1_{10}$ on $C$ which does not satisfy the Gieseker-Petri Theorem and whose residual induces a map $4:1$ on an elliptic curve.\vspace{0.3cm}

\textbf{(IV):} $\deg\phi_{K_C\otimes L^{-1}}=2$. The degree of $\Gamma$ is $6$ and by the Riemann-Hurwitz Formula it follows that $g(\Gamma)\leq 6$. We can assume that $\Gamma$ has only double points as singularities because otherwise $\Gamma$ has a $g^1_k$ for some $k\leq 3$ and, by Remark \ref{rem:secondo}, the curve $[C]\in GP^1_{12,7}$. If $g(\Gamma)\leq 4$, it is easy to check that $\Gamma$ has always a $g^1_3$ and so $[C]\in GP^1_{12,7}$. As a consequence the only two cases that require a more detailed analysis are $g(\Gamma)=5$ and $g(\Gamma)=6$.

Let us consider the case when $\Gamma$ is a plane sextic of genus $5$. We can assume that the singularities of $\Gamma$ are $5$ double points $P_1,\ldots ,P_5$. Some of the $P_i$'s may coincide; indeed, if we need $k$ blow-ups in order to desingularize $\Gamma$ in $P_i$, then this point appears $k$ times in the list. We denote by $x_i, y_i$ the counterimage of $P_i$ under the normalization map $p:Y\rightarrow \Gamma$. Denoting by $B$ and $R$ the branch locus and the ramification locus respectively, the Riemann-Hurwitz Formula implies that both $B$ and $R$ have degree $6$. The double covering $f:C\rightarrow Y$ induced by $\phi$ is given by means of a divisor $\eta$ on $Y$ of degree $-3$ such that $2\eta=-B$ and $f_*\mathcal{O}_C=\mathcal{O}_Y\oplus \mathcal{O}_Y(\eta)$. As $\Pic^{-3}(Y)=Y-Y_4$, we can write $\eta=x-D_4$. 

We consider the divisor $f^*(D_4)\in \Pic^8(C)$. We can assume that $$h^0(C,\mathcal{O}_C(f^*D_4)))=h^0(Y,\mathcal{O}_Y(D_4))+h^0(Y,\mathcal{O}_Y(D_4+\eta))=2,$$ because otherwise we can conclude that $[C]\in \M^2_{12,8}\subset GP^1_{12,7}$. We would like to prove that $\ker\mu_{0,\mathcal{O}_C(f^*D_4)}\neq 0$, which implies $[C]\in GP^1_{12,8}$. By the Base Point Free Pencil Trick we know that $\ker\mu_{0,\mathcal{O}_C(f^*D_4)}\cong H^0(C,K_C\otimes \mathcal{O}_C(f^*D_4)^{-2})$, and this has dimension equal to
$$
h^0(C,f^*(K_Y\otimes\mathcal{O}_Y(-\eta-2D_4)))=h^0(Y,K_Y\otimes\mathcal{O}_Y(-\eta-2D_4))+h^0(Y,K_Y\otimes\mathcal{O}_Y(-2D_4));
$$
here we have used that $K_C=f^*(K_Y\otimes\mathcal{O}_Y(-\eta))$.\\
Since $h^0(Y,K_Y\otimes\mathcal{O}_Y(-2D_4))\neq 0$ whenever $D_4$ is a theta characteristic on $Y$, our goal is to show that $h^0(Y,K_Y\otimes\mathcal{O}_Y(-\eta-2D_4))>0$. As
$$
K_Y\otimes\mathcal{O}_Y(-\eta-2D_4)= \mathcal{O}_Y(3)(-x_1-y_1-\ldots-x_5-y_5-D_4-x),
$$
we need to prove the existence of a plane cubic passing through the points $P_1,P_2$, $P_3$, $P_4$, $P_5$, $p(x)$, $p(z_1)$, $p(z_2)$, $p(z_3)$, $p(z_4)$, where $D_4=z_1+\ldots+z_4$.\\
We can assume that every $g^2_6$ on $Y$ is base point free and not composed with an involution and that every plane model of $Y$ as a sextic has only double points as singularities (otherwise $Y$ would have a $g^1_3$ and $C$ a $g^1_6$); the same is true for all the curves in a neighborhood $U$ of $Y$ in $\M_5$. Up to shrinking $U$, we can assume the existence of a proper morphism $\xi:\mathcal{G}^2_6\rightarrow U$ , where $\mathcal{G}^2_6$ parametrizes couples $(C,l)$, with $[C]\in U$ and $l$ a $g^2_6$ on $C$. We denote by $V_{5,6}$ the variety of irreducible plane curves of degree $6$ and genus $5$ and by $m:V_{5,6}\rightarrow \M_5$ the natural morphism. The locus $m^{-1}(U)$ can be seen as a $PGL(2)$-bundle on $\mathcal{G}^2_6$ parametrizing couples $((C,l),\mathcal{B})$ with $(C,l)\in\mathcal{G}^2_6$ and $\mathcal{B}$ a frame of $l$. Indeed, giving $l$ and $\mathcal{B}$ is equivalent to fixing a plane model of $C$. We denote by $p_1:m^{-1}(U)\rightarrow \mathcal{G}^2_6$ the natural morphism. The restriction $m_U:m^{-1}(U)\rightarrow U$ coincides with the composition $\xi\circ p_1$ and it is proper because both $\xi$ and $p_1$ are. Denoting by $\pi:\M_{5,5}\rightarrow \M_5$ the forgetful map, the morphism 
$$
m_1:m^{-1}(U)\times_U\pi^{-1}(U)\rightarrow\pi^{-1}(U)
$$
is proper because of the invariance of properness under base extension. A point of $m^{-1}(U)\times_U\pi^{-1}(U)$ is of the form $(\Gamma,(C,z_1,\ldots z_5))$, where $C$ is the normalization of $\Gamma$.

We remark that $m^{-1}(U)\times_U\pi^{-1}(U)$ has dimension equal to 
$$
\dim\,\pi^{-1}(U)+\rho(5,2,6)+\dim\,PGL(2)=\dim\,\pi^{-1}(U)+10.
$$
Let $$\mathcal{E}:=H^0(\mathcal{O}_{\mathbb{P}^2}(3))\times (m^{-1}(U)\times_U\pi^{-1}(U))$$ be the trivial bundle on $m^{-1}(U)\times_U\pi^{-1}(U)$ and let us define $\mathcal{F}$ to be the bundle on $m^{-1}(U)\times_U\pi^{-1}(U)$ with fiber over $(\Gamma,(C,z_1\ldots,z_5))$ being the space
$$
H^0(\mathcal{O}_{\mathbb{P}^2}(3)\otimes\mathcal{O}_{\Delta_\Gamma})\oplus\bigoplus_{i=1}^5H^0(\mathcal{O}_{\mathbb{P}^2}(3)\otimes \mathcal{O}_{\phi(z_i)}),
$$
where $\Delta_\Gamma$ is the scheme of all singular points of $\Gamma$. If $\Gamma$ is generic this space is
$$
H^0(\mathcal{O}_{\mathbb{P}^2}(3)\otimes\mathcal{O}_{P_1})\oplus \ldots\oplus H^0(\mathcal{O}_{\mathbb{P}^2}(3)\otimes\mathcal{O}_{P_5})\oplus \bigoplus_{i=1}^5 H^0(\mathcal{O}_{\mathbb{P}^2}(3)\otimes\mathcal{O}_{\phi(z_i)}),
$$
where $P_1,\ldots,P_5$ are the nodes of $\Gamma$. We consider the evaluation map $F:\mathcal{E}\rightarrow \mathcal{F}$. Both $\mathcal{E}$ and $\mathcal{F}$ have rank $10$ and so the degeneracy locus $X(F)$, if nonempty, has codimension $1$ in $m^{-1}(U)\times_U\pi^{-1}(U)$. 

In order to show that $X(F)\neq \emptyset$ we observe that, given a cubic $\Gamma_3\subset\mathbb{P}^2$ and $P_1,\ldots,P_{10}$ ten points on it, one can always find a sextic $\Gamma_6\subset\mathbb{P}^2$ passing through $P_6,\ldots,P_{10}$ and having nodes in $P_1\ldots,P_5$ (because there exists a $\mathbb{P}^{27}$ of plane sextics). Denoting by $\tilde{\phi}:\tilde{C}\rightarrow\Gamma_6$ the normalization map, the point $(\Gamma_6,(\tilde{C},\tilde{\phi}^*(P_6),\ldots,\tilde{\phi}^*(P_{10})))$ lies in $X(F)$. Thus we have that
$$
\dim\, X(F)=\dim\, m^{-1}(U)\times_U\pi^{-1}(U)-1=\dim\,\pi^{-1}(U)+9.$$
As $m_1$ is proper, it follows that $m_1(X(F))$ is closed inside $\pi^{-1}(U)$. Moreover,
$$
\dim\,m_1(X(F))=\dim\,X(F)-\dim\, X_e=\dim\, \pi^{-1}(U)+9-\dim\,X_e,$$
where $X_e$ is the generic fiber of $m_1\vert_{X(F)}$. Therefore $\dim\, m_1(X(F))<\dim\, \pi^{-1}(U)$ if and only if $\dim\,X_e=10$, that is, the generic fiber of $m_1\vert_{X(F)}$ coincides with the generic fiber of $m_1$. If we prove that this cannot happen, then $m_1\vert_{X(F)}$ is surjective and in particular $(Y,p(x),p(z_1),\ldots,p(z_4))\in m_1(X(F))$, which implies the existence of a plane model $\widetilde{\Gamma}$ of $Y$ and of a cubic passing through the singular points of $\widetilde{\Gamma}$ and through the images in $\widetilde{\Gamma}$ of $x,z_1,\ldots,z_4$. Therefore it survives only to prove that $\dim\,X_e\neq 10$. 

Given a general $[C]\in U$, we have to find general points $z_1,\ldots,z_5\in C$, a $g^2_6$ on $C$, together with a frame $\mathcal{B}$ corresponding to a rational map $\phi:C\rightarrow \Gamma\subset\mathbb{P}^{2}$, such that $\Gamma$ has $5$ nodes $P_1,\ldots,P_5$ and there does not exist a cubic through $P_1,\ldots,P_5,\phi(z_1),\ldots,\phi(z_5)$. We remark that any complete $g^2_6$ on $C$ is of the form 
$$L=K_C\otimes\mathcal{O}_C(-a-b),\,\,\,a,b\in C.$$
Having chosen a frame for $H^0(C,L)$ and denoted by $\phi:C\rightarrow \Gamma\subset\mathbb{P}^2$ the corresponding morphism, this is equivalent to saying that
$$
\phi^*\mathcal{O}_\Gamma(1)=\phi^*(\mathcal{O}_\Gamma(3)(-\Delta_\Gamma))\otimes\mathcal{O}_C(-a-b),$$
that is, every cubic in $\mathbb{P}^2$ passing through the singular points of $\Gamma$ and the points $\phi(a),\phi(b)$, intersects $\Gamma$ in other points which are collinear. Choose  $\mathcal{B}$ any frame of $K_C(-z_1-z_2)$; it is enough to take $z_3,z_4,z_5$ such that $\phi(z_3),\phi(z_4),\phi(z_5)$ are not collinear in the plane model of $C$ corresponding to $(K_C\otimes\mathcal{O}_C(-z_1-z_2),\mathcal{B})$.\vspace{0.5cm}

Now we consider the case when $g(\Gamma)=6$, namely $\Gamma$ is a plane sextic with $4$ double points $P_1,\ldots,P_4$. Using the notation introduced above, we now have that $B$ has degree $2$ and so $\eta\in \Pic^{-1}(Y)$. Choose a point $P\in Y$. Since $\Pic^{-2}(Y)=Y_2-Y_4$, we can always write $\eta-P=D_2-D_4$; it follows that $\eta=D_3-D_4$ with $P$ a point of $D_3$. As in the previous case, we can assume that
$$h^0(C,\mathcal{O}_C(f^*D_4))=h^0(Y,\mathcal{O}_Y(D_4))+h^0(Y,\mathcal{O}_Y(D_3))=2,$$
and so $f^*(D_4)$ defines a $g^1_8$ on $C$. In trying to prove that it does not satisfy the Gieseker-Petri Theorem, the above method is unsuccessful. Indeed, we should prove the existence of a plane cubic passing through $P_1,\ldots,P_4$, $p( z_1),\ldots,p(z_6)$, $p(P)$, where $D_4=z_1+\ldots,z_4$, $D_3=z_5+z_6+P$. As $P\in Y$ is arbitrarily chosen, actually it would be enough to prove the existence of a cubic through $P_1,\ldots,P_4,z_1,\ldots,z_6$ and this is a divisorial condition in $(\mathbb{P}^2)^{10}$. Since $\rho(6,2,6)=0$, in this case we do not have any degree of freedom in the choice of a $g^2_6$ on $Y$, namely in the choice of $P_1,\ldots,P_4$.

Thus we proceed in a slightly different way. We have that $\rho(6,2,7)=3$ and, given $l$ a base point free $g^2_7$ on $Y$, we can assume that it defines a birational morphism $$\varphi:Y\rightarrow \Lambda\subset\mathbb{P}^2,$$ where $\Lambda$ is a plane septic of genus $6$; indeed, $l$ cannot be composed with any involution. We expect $\Lambda$ to have only nodes as singularities but in this case we cannot exclude the possibility that $\Lambda$ has some triple points. As $Y$ is the normalization of $\Lambda$, we have that
$$
K_Y=\varphi^*(\mathcal{O}_\Lambda(4)(-\Delta_\Lambda))\textrm{   with   }\Delta_\Lambda=\sum_{P\in\mathrm{Sing}\Lambda}(r_P-1)P,
$$
where $r_P$ is the multiplicity of $\Lambda$ in $P$. Of course for $\Lambda$ generic, the singular locus $\Delta_\Lambda$ is the sum of the nine nodes $P_1,\ldots,P_9$ and the condition $\ker\mu_{0,\mathcal{O}_C(f^*D_4)}\neq 0$ is equivalent to the existence of a plane quartic through $P_1,\ldots,P_9$, $\varphi(z_1),\ldots,\varphi(z_6)$. In the non generic case we have a different condition equivalent to $\ker\mu_{0,\mathcal{O}_C(f^*D_4)}\neq 0$ (for instance, when $\Lambda$ has a triple point $Q$ and six double points $P_1,\ldots,P_6$, then we require that the plane quartic has a double point in $Q$ and passes through $P_1,\ldots,P_6$). However, the number of independent conditions imposed on the plane quartics is the same.

As before, we consider a neighborhood $U$ of $Y$ in $\M_6$ such that there exists a proper morphism $\xi:\mathcal{G}^2_7\rightarrow U$, where $\mathcal{G}^2_7$ parametrizes couples $(C,l)$, with $[C]\in U$ and $l$ a $g^2_7$ on $C$. We can assume that, given $[C]\in U$, the generic $g^2_7$ on $C$ is base point free and not composed with an involution but in this case the models of $C$ as a plane septic can have also some triple points. Denoting by $m:V_{6,7}\rightarrow \M_6$ the natural morphism, the restriction $m_U:m^{-1}(U)\rightarrow U$ is proper. If $\pi:\M_{6,6}\rightarrow \M_6$ is the forgetful map, then the induced morphism $m_1:m^{-1}(U)\times_U\pi^{-1}(U)\rightarrow \pi^{-1}(U)$ is proper, too. We have that
$$
\begin{array}{llll}
\dim\,m^{-1}(U)\times_U\pi^{-1}(U)&=&\dim\,\pi^{-1}(U)+\rho(6,2,7)+\dim\,PGL(2)&=\\
&=&\dim\pi^{-1}(U)+11&.
\end{array}
$$
As in the previous case, we define $$\mathcal{E}:=H^0(\mathcal{O}_{\mathbb{P}^2}(4))\times (m^{-1}(U)\times\pi^{-1}(U))$$ and $\mathcal{F}$ being the bundle over $m^{-1}(U)\times\pi^{-1}(U)$ whose fiber over $(\Lambda,(C,z_1,\ldots,z_6))$ is 
$$H^0(\mathcal{O}_{\mathbb{P}^2}(4)\otimes\mathcal{O}_{\Delta_\Lambda})\oplus H^0(\mathcal{O}_{\mathbb{P}^2}(4)\otimes\mathcal{O}_{\varphi(z_1)})\oplus\ldots\oplus H^0(\mathcal{O}_{\mathbb{P}^2}(4)\otimes\mathcal{O}_{\varphi(z_6)}).$$
For $\Lambda\in V_{6,7}$ generic we have that
$$H^0(\mathcal{O}_{\mathbb{P}^2}(4)\otimes\mathcal{O}_{\Delta_\Lambda})=H^0(\mathcal{O}_{\mathbb{P}^2}(4)\otimes\mathcal{O}_{P_1})\oplus\ldots\oplus H^0(\mathcal{O}_{\mathbb{P}^2}(4)\otimes\mathcal{O}_{P_9}),$$
where $P_1,\ldots,P_9$ are the nodes of $\Lambda$. Instead, if for instance $\Lambda$ has one triple point $Q$ and $6$ nodes $P_1,\ldots,P_6$, then the following equality holds:
$$H^0(\mathcal{O}_{\mathbb{P}^2}(4)\otimes\mathcal{O}_{\Delta_\Lambda})=H^0(\mathcal{O}_{\mathbb{P}^2}(4)\otimes\mathcal{O}_{2Q})\oplus \ldots\oplus H^0(\mathcal{O}_{\mathbb{P}^2}(4)\otimes\mathcal{O}_{P_6}).$$
 We define $F:\mathcal{E}\rightarrow\mathcal{F}$ to be the evaluation map. As both $\mathcal{E}$ and $\mathcal{F}$ have rank $15$, the situation is analogous to the one already treated. Therefore, in order to prove that the image under $m_1$ of the degeneracy locus $X(F)$ is the whole $\pi^{-1}(U)$, it is enough to show that the generic fiber $X_e$ of $m_1\vert_{X(F)}$ is nonempty and that it does not coincide with the generic fiber of $m_1$. The fact that $X_e\neq \emptyset$ follows easily by observing that, given $15$ points on a quartic $\Lambda_4\subset\mathbb{P}^2$, there always exists a plane septic $\Lambda_7$ passing through them and having nodes in the first nine. On the other hand, it can be shown that $\dim\,X_e\neq 15$ by proceeding like in the case of genus $5$ because on a curve $C$ of genus $6$ any complete $g^2_7$ is of the form 
$
l=K_C\otimes\mathcal{O}_C(-a-b-c)$, with $a,b,c\in C$.
\end{proof}
Finally, we obtain that:
\begin{cor}
The locus $GP_{12}$ has pure codimension $1$ in $\M_{12}$.
\end{cor} 
\begin{proof}
By the remarks at the beginning of the section we have to study only the component $GP^1_{12,10}$. Given $[C]\in GP^1_{12,10}$ we can assume that $[C]\not\in GP^1_{12,9}\cup GP^2_{12,10}\cup GP^2_{12,11}$, which forces any $l\in G^1_{10}(C)$ for which the Gieseker-Petri Theorem fails to verify condition $(1)$. Theorem \ref{thm:coverings} implies that if $[C]\not\in GP^1_{12,7}\cup GP^1_{12,8}$, then condition $(2)$ is satisfied, too. We can thus apply Proposition \ref{prop:bo}.
\end{proof}


\begin{thebibliography}{ACGH}
\bibitem[AC1]{abete} E. Arbarello, M. Cornalba, {\em A few remarks about the variety of irreducible plane curves of given degree and genus}, Ann. Scient. \'Ec. Norm. Sup. 4 sreie \textbf{16} (1983), 467-488.


\bibitem[AC2]{abete2}
E. Arbarello, M. Cornalba, {\em Footnotes to a Paper of Beniamino Segre}, Math. Ann. \textbf{256} (1981), 341-362.

\bibitem[ACGH]{alberello}
E. Arbarello, M. Cornalba, P. Griffiths, J. Harris, {\em Geometry of algebraic curves}, Grundl. Math. Wiss. \textbf{267}, Springer Verlag, 1985.

\bibitem[BS]{sernesi}
A. Bruno, E. Sernesi, {\em A note on the Petri loci}, \texttt{ArXiv}:1012.0856v2.


\bibitem[Ca1]{casto}
A. Castorena, {\em Curves of genus seven that do not satisfy the Gieseker-Petri Theorem}, Bollettino U.M.I. \textbf{8-B} (2005), 697-706.

\bibitem[Ca2]{casto2}
A. Castorena, {\em A family of plane curves with moduli $3g-4$}, Glasgow Math. J. \textbf{49} (2007), 417-422.

\bibitem[Ca3]{casto3}
A. Castorena, {\em Remarks on the Gieseker-Petri divisor in genus eight}, Rend. Circ. Mat. Palermo \textbf{59} (2010), 143-150.

\bibitem[EH1]{altra}
D. Eisenbud, J. Harris, {\em A simple proof of the Gieseker-Petri theorem on special divisors}, Invent. Math. \textbf{74} (1983), 269-280

\bibitem[EH2]{harri}
D. Eisenbud, J. Harris, {\em The Kodaira dimension of the moduli space of curves of genus $\geq 23$}, Invent. Math. \textbf{90} (1987), 359-387.

\bibitem[EH3]{eisi}
D. Eisenbud, J. Harris, {\em Irreducibility of some families of linear series with Brill-Noether Number $-1$}, Ann. Scient. \'Ec. Norm. Sup. 4e s\'erie \textbf{22} (1989), 33-53.

\bibitem[ELMS]{cliff}
D. Eisenbud, H. Lange, G. Martens, and F.-O. Schreyer, {\em The Clifford dimension of a projective curve}, Compos. Math. \textbf{72} (1989), 173-204. 

\bibitem[F1]{gabi1}
G. Farkas, {\em Gaussian maps, Gieseker-Petri loci and large theta-characteristics}, J. Reine Angew. Math \textbf{581} (2005), 151-173.

\bibitem[F2]{gabi}
G. Farkas, {\em Rational maps between moduli space of curves and Gieseker-Petri divisor}, J. Algebraic Geom. \textbf{19} (2010), 243-284.

\bibitem[Fu]{fulton}
W. Fulton, {\em Hurwitz schemes and irreducibility of moduli of algebraic curves}, Ann. Math. \textbf{90} (1969), 542-575. 

\bibitem[Gi]{gies}
D. Gieseker, {\em Stable curves and special divisors}, Invent. Math. \textbf{66} (1982), 251-275.

\bibitem[Ha1]{harti}
R. Hartshorne, {\em Ample vector bundles}, Inst. Hautes \'Etudes Sci. Publ. Math. \textbf{36} (1969), 75-110.

\bibitem[Ha2]{harti2}
R. Hartshorne, {\em Algebraic geometry}, Springer-Verlag, New York-Heidelberg-Berlin, 1977.

\bibitem[Ma]{marti}
G. Martens, {\em \"Uber den Clifford-Index algebraischer Kurven}, J. Reine Angew. Math. \textbf{36} (1982), 83-90. 

\bibitem[St]{steffen}
F. Steffen, {\em A generalized principal ideal theorem with an application to Brill-Noether theory}, Invent. Math. \textbf{132} (1998), 73-89.

\bibitem[Te]{tex}
M. Teixidor, {\em The divisor of curves with a vanishing theta-null}, Compos. Math. \textbf{66} (1988), 15-22.
\end{thebibliography}
\end{document}